\documentclass[a4paper,10pt,reqno]{amsart}
\usepackage[centertags]{amsmath}
\usepackage{amsfonts,amssymb,amsthm} 
\usepackage[dvips]{graphicx} 
\usepackage{psfrag}
\usepackage[english]{babel}                
\usepackage{newlfont}
\usepackage[body={15cm, 23.1cm},centering]{geometry} 
\usepackage{fancyhdr}
 \theoremstyle{plain}
 \begingroup
 \newtheorem{teo}{Theorem}[section]
 \newtheorem{lem}[teo]{Lemma}
 \newtheorem{prop}[teo]{Proposition}
 
 \endgroup 
  
 \theoremstyle{definition}
 \begingroup
 \newtheorem{defi}[teo]{Definition}
 
 \endgroup 
  
 \theoremstyle{remark}
 \begingroup
 \newtheorem{rem}[teo]{Remark}
 \endgroup
 
 \renewcommand{\qedsymbol}{$\blacklozenge$}
 \numberwithin{equation}{section}

\renewcommand{\a}{\mathcal{A}}
\renewcommand{\b}{\mathcal{B}}
\newcommand{\f}{\mathcal{L}}
\newcommand{\m}{\mathcal{M}} 
\newcommand{\p}{\mathcal{P}} 
\newcommand{\W}{\mathcal{W}}
\newcommand{\V}{\mathcal{V}}
\newcommand{\D}{\mathcal{D}}
\newcommand{\Y}{\mathcal{Y}}

\newcommand{\e}{\varepsilon}
\newcommand{\Om}{\Omega}

\newcommand{\R}{{\mathbb R}}
\newcommand{\NP}{\mathbb{N}^{+}}

\newcommand{\M}{\mathbb{M}}

\newcommand{\g}{\nabla}
\newcommand{\dig}{\mathrm{diag}}

\renewcommand{\neg}{\negthinspace}
\newcommand{\nek}{\negmedspace}
\newcommand{\nec}{\negthickspace}

\newcommand{\norm}[1]{\left\Vert#1\right\Vert}
\newcommand{\snorm}[1]{\left|#1\right|}

\begin{document}
\title[Loss of polyconvexity by homogenization]
{Loss of polyconvexity by homogenization: a new example}
\author[Marco Barchiesi]{Marco Barchiesi}
\date{30th March 2006}
\address[Marco Barchiesi]{S.I.S.S.A., Via Beirut 2-4, 34014, Trieste, Italy}
\email[Marco Barchiesi]{barchies@sissa.it}
\maketitle
%%----------------------------------------------------------------------------------------------------------------------

\begin{center}
\begin{minipage}{12cm}
\small{
\noindent {\bf Abstract.} 
This article is devoted to the study of the asymptotic behavior of the zero-energy deformations set of a 
periodic nonlinear composite material. We approach the problem using two-scale Young measures. 
We apply our analysis to show that polyconvex energies are not closed with respect to periodic
homogenization. The counterexample is obtained through a rank-one laminated structure 
assembled by mixing two polyconvex functions with $p$-growth, where $p\geq2$ can be fixed arbitrarily.

\vspace{10pt}
\noindent {\bf Keywords:}  
composite materials, homogenization, quasiconvexity, polyconvexity,
rank-one laminates, two-scale Young measures, $\Gamma$-convergence

\vspace{6pt}
\noindent {\bf 2000 Mathematics Subject Classification:} 28A20, 35B27, 35B40, 49J45, 73B27}
%28A20 funzioni misurabili e modi di convergenza
%35B27 omogenizzazione con struttura periodica
%35B40 comportamento asintotico delle soluzioni
%49J45 metodi di semicontinuita` e convergenza
%73B27 materiali non omogenei, omogenizzazione
\end{minipage}
\end{center}

\bigskip
\tableofcontents
%%----------------------------------------------------------------------------------------------------------------------

\section{Introduction}

\noindent
Many problems related to composite materials (see \cite{Mil02} and references therein)
lead to the variational analysis of families of integral functionals.

In the case of a periodic composite, the energy density can be described by a family of the type
$\lbrace(\W_k,P_k)\rbrace_{k=1,\ldots,n}$, where
$\lbrace\W_k:\M^{s\times d}\nek\rightarrow\nek[0,+\infty)\rbrace_{k=1,\ldots,n}$ is a family of continuous functions describing the energy densities of the components and $\lbrace P_k\rbrace_{k=1,\ldots,n}$ is a 
measurable partition of the unit cell $\square:=[0, 1)^d$.
In a mix of fineness~$\e$, the functional $E_\e$ that represents, at microscopic level, 
the stored energy of the composite is of the type 
\begin{equation*}
E_\e (u)=\int_\Om \W\Bigl(\langle\frac{x}{\e}\rangle,\g u(x)\Bigr)dx,
\end{equation*}
where $\langle\cdot\rangle$ denotes the fractional part of a vector componentwise,
$\Om$ is an open bounded domain in $\R^d$ and $\W(y,\Lambda)=\sum_{k=1}^{n} \chi_{P_k}(y)\,\W_k(\Lambda)$
with $\chi_{P_k}$ the characteristic function of $P_k$. 

When the parameter $\e$ tends to $0$, the microscopic structure becomes finer and finer
and the asymptotic behaviour of the composite is  that of a homogeneous material.
If the function $\W$ satisfies coerciveness and growth conditions
of order $p\in(1,+\infty)$, \textit{i.e.}, there exist $c_1$, $c_2$, $c_3>0$ such that
\begin{equation}\label{crescita}
c_1\snorm{\Lambda}^p-c_2\leq\W\left(y,\Lambda\right)\leq c_3\left(1+\snorm{\Lambda}^p\right)
\quad\text{for all} \;(y,\Lambda)\in\square\times\M^{s\times d}\,,
\end{equation}
then the stored energy $E_{hom}$ of this material can be described efficiently 
by the $\Gamma$-limit of the family $E_\e$ and is of the type
\begin{equation*}
E_{hom} (u)=\int_\Om \W_{hom}\bigl(\g u(x)\bigr)dx,
\end{equation*}
where $\W_{hom}$ is a suitable quasiconvex function called \emph{homogenized integrand}. 

\vspace{5pt}
\noindent
\textbf{Question}:
Suppose that the sets $\mathfrak{A}_k:=\W_k^{-1}\left(0\right)$ ($k=1,\ldots,n$) are not empty.
Is it possible to characterise $\W_{hom}^{-1}\left(0\right)$ and its dependence on
$\lbrace(\mathfrak{A}_k,P_k)\rbrace_{k=1,\ldots,n}$?

\noindent
\begin{defi}
Given a measurable partition $\lbrace P_k\rbrace_{k=1,\ldots,n}$ of the unit cell $\square$
and a family of compact sets $\lbrace\mathfrak{A}_k\rbrace_{k=1,\ldots,n}$,
we define $\mathfrak{A}_{hom}$ as the set of matrices $A\in\M^{s\times d}$ such that
for all $\e_h\rightarrow0^+$ there exists a sequence $u_h\in W^{1,\infty}(\Omega,\R^s)$ satisfying
\begin{itemize}
\item[i)]   $u_h\rightharpoonup A\,x$ \;weakly* in $W^{1,\infty}(\Omega,\R^s)$;
\item[ii)]  $\sum_{k=1}^{n}\chi_{P_k}\bigl(\langle\frac{x}{\e_h}\rangle\bigr)
       \mathrm{dist}\bigl(\g u_h(x),\mathfrak{A}_k\bigr)\rightarrow 0$ in measure.
\end{itemize}
The second condition is equivalent (see proposition \ref{local supporti}) to the statement that the 
two-scale gradient Young measure $(\nu_{(x,y)})_{(x,y)\in\Om\times\square}$ corresponding to 
(a subsequence of) $\nabla u_h$ satisfies 
\begin{itemize}
\item[ii)'] for all \;$k\in\lbrace1,\ldots,n\rbrace$ and for \textit{a.e.} $(x,y)\in\Om\times P_k$
we have \;$\mathrm{supp}\,\nu_{(x,y)}\subseteq\mathfrak{A}_k$ .
\end{itemize}
We call $\mathfrak{A}_{hom}$ the \emph{homogenized set} related to $\lbrace(\mathfrak{A}_k,P_k)\rbrace_{k=1,\ldots,n}$.
\end{defi}

\vspace{4pt}
The first result of this article gives an answer to the previous question, 
showing the link between homogenized integrands and homogenized sets.
We were inspired by \cite[section 4]{B92} and \cite[section~4]{B02}.
\begin{teo}\label{zero energy}
Let $\lbrace\W_k:\M^{s\times d}\nek\rightarrow\nek[0,+\infty)\rbrace_{k=1,\ldots,n}$ be 
a family of continuous functions, let $\lbrace P_k\rbrace_{k=1,\ldots,n}$ be a measurable 
partition of the unit cell $\square$ and let $\W(y,\Lambda):=\sum_{k=1}^{n} \chi_{P_k}(y)\,\W_k(\Lambda)$.
Suppose that for $k\in\lbrace1,\ldots,n\rbrace$ 
\begin{itemize}
\item[i)] the function $\W_k$ satisfies coerciveness and growth conditions of order $p$;
\item[ii)] the set $\mathfrak{A}_k:=\W_k^{-1}(0)$ is not empty.
\end{itemize}
Then
\vspace{-2pt}
\begin{equation*}
\mathfrak{A}_{hom}=\W_{hom}^{-1}(0),
\end{equation*}

\vspace{3pt}
\noindent
\textit{i.e.}, the zero-level set of the homogenized integrand $\W_{hom}$ related to $\W$ coincide with
the homogenized set $\mathfrak{A}_{hom}$ related to $\lbrace(\mathfrak{A}_k,P_k)\rbrace_{k=1,\ldots,n}$.
\end{teo}
Let observe that in the simple case $n=1$, the function $\W_{hom}$ is the quasiconvexification 
of $\W_1$ and $\mathfrak{A}_{hom}$ is the quasiconvex hull of $\mathfrak{A}_1$.

\vspace{6pt}
The second result is to show that the polyconvexity of the sets $\mathfrak{A}_k$ ($k=1,\ldots,n$)
does not ensure the polyconvexity of $\mathfrak{A}_{hom}$. 
This aim is achieved with an explicit example in the case $d=s=n=2$.
\begin{teo}\label{main}
\textnormal{(Loss of polyconvexity)}. Consider the sets
\begin{itemize}
\item $P_1=[0,\frac{1}{2})\times[0,1)$ and $P_2=[\frac{1}{2},1)\times[0,1)$;
\item $\mathfrak{A}_1=\lbrace O, A^1_1, A^2_1\rbrace$ 
      with $O=\dig(0,0)$, $A^1_1=\dig(-\frac{1}{2},1)$ and $A^2_1=\dig(-2,2)$;
\item $\mathfrak{A}_2=\lbrace O, A^1_2, A^2_2\rbrace$ 
      with $A^1_2=\dig(\frac{5}{2},1)$ and $A^2_2=\dig(3,2)$;
\item $\mathfrak{B}=\lbrace O, B^1, B^2\rbrace$ with $2B^j=A^j_1+A^j_2$ for $j=1,2$.
\end{itemize}
The following properties hold:
\begin{itemize}
\item[i)] $\mathfrak{A}_1$ and $\mathfrak{A}_2$ are both polyconvex sets;
\item[ii)] fixed $p\geq2$, there exist two polyconvex functions
           \,$\W_1,\W_2:\M^{2\times2}\nek\rightarrow\nek[0,+\infty)$ $p$-coercive and with $p$-growth such that 
           $\mathfrak{A}_1=\W_1^{-1}(0)$ and $\mathfrak{A}_2=\W_2^{-1}(0)$;
\item[iii)] the homogenized set \;$\mathfrak{A}_{hom}$ related to $\lbrace(\mathfrak{A}_k,P_k)\rbrace_{k=1,2}$ 
            is not polyconvex. More precisely
\begin{equation*}
\mathfrak{B}\subseteq\mathfrak{A}_{hom} \quad\text{and}\quad\; \mathfrak{B}^{pc}\nsubseteq\mathfrak{A}_{hom},
\end{equation*} 
where $\mathfrak{B}^{pc}$ denotes the polyconvex hull of $\mathfrak{B}$;
\item[iv)] the homogenized integrand $\W_{hom}$ related to $\W(y,\Lambda):=\chi_{P_1}(y)\,\W_1(\Lambda)+\chi_{P_2}(y)\,\W_2(\Lambda)$
           is not polyconvex.
\end{itemize} 
\end{teo}  

\vspace{5pt}
\noindent
The example given in the result above shows that polyconvexity is not preserved by homogenization,
unlike convexity and quasiconvexity. 

The first example of loss of polyconvexity by homogenization is due to Braides \cite{Br94}.
He considers a function $\W$ assembled by using two functions 
\,$\W_1,\W_2:\M^{2\times2}\nek\rightarrow\nek[0,+\infty)$ with different growth conditions.
More precisely
\begin{equation*}
\W(y,\Lambda):=
\begin{cases}
\W_1(\Lambda) & \text{if} \quad y\in\square\backslash[\frac{1}{4},\frac{3}{4}]^2,\\
\W_2(\Lambda) & \text{if} \quad y\in[\frac{1}{4},\frac{3}{4}]^2,
\end{cases}
\end{equation*}
where $\W_1$ is a convex function satisfying coerciveness and growth conditions of order $p<2$, 
while $\W_2$ is a polyconvex function satisfying coerciveness and growth conditions of order $2$.
Since the homogenized integrand related to $\W$ is not convex and satisfies a growth condition of order
$p<2$, it cannot be polyconvex. A suitable quadratic perturbation, considered in \cite{Br98}, permits
to assume that also $\W_1$ satisfies coerciveness and growth conditions of order $2$.

Our example is quite different. It is based only on the structure of $\lbrace(\mathfrak{A}_k,P_k)\rbrace_{k=1,2}$
and is independent of the growth condition.
Indeed, we first construct two sets $\mathfrak{A}_1$ and $\mathfrak{A}_2$ such that the homogenized set
\;$\mathfrak{A}_{hom}$ related to $\lbrace(\mathfrak{A}_k,P_k)\rbrace_{k=1,2}$ is not polyconvex
and then, fixed $p\neg\geq\neg2$, we construct two polyconvex functions
\,$\W_1,\W_2:\M^{2\times2}\nek\rightarrow\nek[0,+\infty)$ 
$p$-coercive and with $p$-growth such that $\mathfrak{A}_1=\W_1^{-1}(0)$ and $\mathfrak{A}_2=\W_2^{-1}(0)$.
The homogenized integrand $\W_{hom}$ related to
$\W(y,\Lambda):=\chi_{P_1}(y)\,\W_1(\Lambda)+\chi_{P_2}(y)\,\W_2(\Lambda)$
cannot be polyconvex because $\mathfrak{A}_{hom}=\W_{hom}^{-1}(0)$.

The idea behind our construction of $\lbrace(\mathfrak{A}_k,P_k)\rbrace_{k=1,2}$ is the following.
The set $\mathfrak{A}_1$ (resp. $\mathfrak{A}_2$) is polyconvex because the difference of two elements
of $\mathfrak{A}_1$ (resp. $\mathfrak{A}_2$) has negative (resp. positive) determinant.
The set $\mathfrak{B}=\lbrace O,I,\dig(\frac{1}{2},2)\rbrace$, composed of average of the
correspondent matrices of $\mathfrak{A}_1$ and $\mathfrak{A}_2$, does not have these properties
and its polyconvex hull is not trivial.
The fact that $\mathfrak{B}^{pc}\nsubseteq\mathfrak{A}_{hom}$ is proved by using the function $\V$
defined in \eqref{funz sverak}. This function was used originally by \v Sver\'ak in \cite{Sver92}
to prove the quasiconvexity of sets of the type  $\lbrace O,I,\dig(b_1,b_2)\rbrace$, where
$0<b_1<1$ and $b_2>1$. See also \cite{AsFa02} for a different proof.

\vspace{4pt}
The plan of the paper is as follows.
In section 2 we collect concepts and basic facts about Young measures and two-scale Young measures.
In section 3 we recall some well known facts about $\Gamma$-convergence and prove Theorem \ref{zero energy}.
In section 4 we present some simple results about polyconvexity.
Finally, in section 5 we provide a proof of Theorem \ref{main}.
%%%%%%%%%%%%%%%%%%%%%%%%%%%%%%%%%%%%%%%%%%%%%%%%%%%%%%%%%%%%%%%%%%%%%%%%%%%%%%%%%%%%%%%%%%%%%%%%%%%
%%%%%%%%%%%%%%%%%%%%%%%%%%%%%%%%%%%%%%%%%%%%%%%%%%%%%%%%%%%%%%%%%%%%%%%%%%%%%%%%%%%%%%%%%%%%%%%%%%%
\section{Two-scale Young measures}
\noindent
We start by collecting preliminary results about Young measures (see \cite{B88}, \cite{Mul99bis} and 
\cite{Val90}) and two-scale Young measures (see \cite{Bar05} and \cite{Val97}).

Before we list the notation used in the following:
\begin{itemize}
\item $\Omega$ a bounded open subset of $\R^d$;
\item $\square$ the unit cell $[0, 1)^d$; 
\item $\f(\Om)$ the Lebesgue $\sigma$-algebra on $\Om$\,;
\item $\snorm{S}$ the Lebesgue measure of a measurable subset $S\subseteq\R^l$;
\item $\b(S)$ the Borel $\sigma$-algebra on a subset $S\subseteq\R^l$;
\item $\p(\R^m)$ the set of probability measures on $\R^m$;
\item $\Y(\Om,\R^m)$ the family of all weakly* measurable maps 
$x\in\Om\negthinspace\xrightarrow{\mu}\negthinspace\mu_x\in\p(\R^m)$;
\\a corresponding definition holds for $\Y(\Om\times\square,\R^m)$;
\item $\langle x\rangle\in\square$ the fractional part of $x\in\R^d$ componentwise, \textit{i.e.}, 
\begin{equation*}
\langle x\rangle_k=x_k-\lfloor x_k\rfloor \quad\text{for} \;k\in\lbrace1,\ldots,d\rbrace,
\end{equation*}
where $\lfloor x_k\rfloor$ stands for the largest integer less than or equal to $x_k$;
\item $\e_h$ a sequence in $(0,+\infty)$.
 \end{itemize}

\begin{rem}
\begin{itemize}\label{oss princ}
\item[]
\item[i)] The term \emph{measurable} is tacitly understood as Lebesgue measurable, unless otherwise stated.
\item[ii)]A map $\mu:\Om\rightarrow\p(\R^m)$ is said to be \emph{weakly* measurable} if 
         $x\rightarrow\mu_x(S)$ is measurable for all $S\in\b(\R^m)$.
         By approximation, if $\mu\in\Y(\Om,\R^m)$ and $\W:\Om\times\R^m\rightarrow[0,+\infty)$
         is $\f(\Om)\otimes\b(\R^m)$-measurable, then $x\rightarrow\int_{\R^m}\W(x,\lambda)d\mu_x(\lambda)$
         is measurable.
\vspace{2pt}
\item[iii)]More precisely, the elements of $\Y(\Om,\R^m)$ are equivalence
     classes of maps that agree \textit{a.e.}; we usually do not distinguish 
     these maps from their equivalence classes.
\end{itemize}
\end{rem} 
 
\vspace{6pt}
The following result is known as the Fundamental Theorem on Young Measures.
It shows that the weak limit of a sequence of the type $\W\left(\cdot, u_{h}(\cdot)\right)$
can be expressed through a suitable map $\mu\in\Y(\Om,\R^m)$ associated to $u_h$.
The proof can be found in \cite{B88} (see also \cite{Mul99bis} and \cite{Val90}).
We recall that a $\f(\Om)\otimes\b(\R^m)$-measurable function $\W$ is a Carath\'eodory integrand if
$\W(x,\cdot)$ is continuous for \textit{a.e.} $x\in\Om$. 

\begin{teo}\label{fund young}
Let $u_h$ be a bounded sequence in $L^1(\Om,\R^m)$. There exist a subsequence $u_{h_i}$ and a 
map $\mu\in\Y(\Om,\R^m)$ such that the following properties hold:
\begin{enumerate}
\item[i)] if $\W:\Om\times\R^m\rightarrow[0,+\infty)$ is a Carath\'eodory integrand, then
\begin{equation*}
\liminf_{i\rightarrow+\infty}\int_\Om\W\bigl(x, u_{h_i}(x)\bigr)dx\geq\int_\Om\overline{\W}(x)dx
\end{equation*}
where
\begin{equation*}
\overline{\W}(x):=\int_{\R^m}\W(x,\lambda)d\mu_x(\lambda);
\end{equation*}
\item[ii)] if $\W:\Om\times\R^m\rightarrow\R$ is a Caratheodory integrand and
$\W\left(\cdot, u_{h}(\cdot)\right)$ is equi-integrable, then
$\W(x,\cdot)$ is $\mu_x$-integrable for \textit{a.e.} $x\in\Om$, \;$\overline{\W}$ is in $L^1(\Om)$ and
\begin{equation*}
\lim_{i\rightarrow+\infty}\int_\Om\W\bigl(x, u_{h_i}(x)\bigr)dx=\int_\Om\overline{\W}(x)dx;
\end{equation*}
\item[iii)] if $\mathfrak{A}\subseteq\R^m$ is closed, then $\mathrm{supp}\,\mu_x\subseteq\mathfrak{A}$ for \textit{a.e.} $x\in\Om$ if and only if $\mathrm{dist}(u_{h_i},\mathfrak{A})\rightarrow 0$ in measure.
\end{enumerate}
\end{teo}
\begin{defi}
The map $\mu$ is called the \emph{Young measure} generated by the sequence $u_{h_i}$.
\end{defi}
%%%%%%%%%%%%%%%%%%%%%%%%%%%%%%%%%%%%%%%%%%%%%%%%%%%%%%%%%%%%%%%%%%%%%%%%%%

\vspace{6pt}
In homogenization processes, we are interested to asymptotic behaviour
of sequences of the type $\W_h(\cdot)=\W\bigl(\cdot,\langle\frac{\cdot}{\e_h}\rangle, u_h(\cdot)\bigr)$.
Under technical hypothesis on $\W$, the weak limit of $\W_h$ can be expressed through a suitable map 
$\nu\in\Y(\Om\times\square,\R^m)$ associated to $u_h$. 

\begin{defi}\label{admiss}
A function $\W:\Om\times\square\times\R^m\rightarrow\R$ 
is said to be an \emph{admissible integrand}
if there exist a family $\lbrace X_\delta\rbrace_{\delta>0}$ of compact subsets of $\Om$ 
and a family $\lbrace Y_\delta\rbrace_{\delta>0}$ of compact subsets of $\square$ such that 
$\snorm{\Om\backslash X_\delta}\leq\delta$, $\snorm{\square\backslash Y_\delta}\leq\delta$ 
and $\W\vert _{X_\delta\times Y_\delta\times\R^m}$ is continuous for every $\delta>0$. 
\end{defi}

\begin{rem}
It is not restrictive to suppose that the families $\lbrace X_\delta\rbrace_{\delta>0}$
and $\lbrace Y_\delta\rbrace_{\delta>0}$ are decreasing, \textit{i.e.}, $\delta'\leq\delta$ 
implies $X_{\delta}\subseteq X_{\delta'}$ and $Y_{\delta}\subseteq Y_{\delta'}$. 
Otherwise, it is sufficient to consider the new families $\lbrace\widetilde{X}_\delta\rbrace_{\delta>0}$
and $\lbrace\widetilde{Y}_\delta\rbrace_{\delta>0}$, where
\begin{equation*}
\widetilde{X}_{\delta}:=\bigcap_{i\geq i_{\delta}}X_{2^{-i}}\,, 
\quad \widetilde{Y}_{\delta}:=\bigcap_{i\geq i_{\delta}}Y_{2^{-i}}
\end{equation*}
and $i_{\delta}$ is the minimum positive integer such that $2^{1-i_{\delta}}\leq\delta$.
\end{rem}
\begin{rem}\label{misurab}
Admissible integrands have good measurability properties:
if $\W$ is an admissible integrand, then
there exist a Borel set $X\subseteq\Om$ with $\snorm{\Om\backslash X}=0$ 
and a Borel set $Y\subseteq\square$ with $\snorm{\square\backslash Y}=0$, such that
$\W\vert _{X\times Y\times\R^m}$ is borelian.
In particular, for every fixed $\e$, the function $(x,\lambda)\rightarrow\W\left(x,\langle\frac{x}{\e}\rangle,\lambda\right)$
is $\f(\Om)\otimes\b(\R^m)$-measurable. 
\end{rem}
\begin{rem}\label{scorza}
Let $\lbrace\W_k:\Om\times\M^{s\times d}\nek\rightarrow\nek[0,+\infty)\rbrace_{k=1,\ldots,n}$ be 
a family of Carath\'eodory integrands and let $\lbrace P_k\rbrace_{k=1,\ldots,n}$ be a measurable 
partition of the unit cell $\square$. By applying Lusin theorem to each $\chi_{P_k}$ and 
Scorza-Dragoni theorem (see \cite{Eke99}) to each $\W_k$, we obtain that 
$\W(x,y,\lambda):=\sum_{k=1}^{n} \chi_{P_k}(y)\,\W_k(x,\lambda)$ is an admissible integrand.
\end{rem}

\vspace{3pt}
The next two results are the equivalent of Theorem \ref{fund young} in the case of two-scale Young measures.
\begin{teo}\label{fund multi young}
Let $\e_h\rightarrow0^+$, let $u_h$ be a bounded sequence in $L^1(\Om,\R^m)$ and 
let $w_h:\Om\rightarrow\square\times\R^m$ be defined by $w_h(x):=\bigl(\langle\frac{x}{\e_h}\rangle, u_h(x)\bigr)$.
There exist a subsequence $w_{h_i}$ and a map $\nu\in\Y(\Om\times\square,\R^m)$ such that 
the following properties hold:
\begin{enumerate}
\item[i)] if $\W:\Om\times\square\times\R^m\rightarrow[0,+\infty)$ is an admissible integrand, then
\begin{equation}\label{semicont}
\liminf_{i\rightarrow+\infty}\int_{\Om}
\W\bigl(x, w_{h_i}(x)\bigr)dx\geq\int_{\Om\times\square}\overline{\W}(x,y)\,dx\,dy ,
\end{equation}
where 
\begin{equation*}
\overline{\W}(x,y):=\int_{\R^m}\W(x,y,\lambda)\;d\nu_{(x,y)}(\lambda);
\end{equation*}
\item[ii)] if $\W:\Om\times\square\times\R^m\rightarrow\R$ is an admissible integrand and
$\W(\cdot,w_{h}(\cdot))$ is equi-integrable, then $\W(x,y,\cdot)$ is $\nu_{(x,y)}$-integrable 
for \textit{a.e.} $(x,y)\in\Om\times\square$, $\overline{\W}$ is in $L^1(\Om\times\square)$ and
\begin{equation}\label{cont}
\lim_{i\rightarrow+\infty}\int_{\Om}
\W\bigl(x, w_{h_i}(x)\bigr)dx=\int_{\Om\times\square}\overline{\W}(x,y)\,dx\,dy.
\end{equation}
\end{enumerate}
\end{teo}
\begin{defi}
The map $\nu$ is called the \emph{two-scale Young measure} generated by the sequence $u_{h_i}$
with respect to $\e_{h_i}$. In the sequel we omit to specify the dependence on $\e_{h_i}$ 
when it is clear from the context.
\end{defi}
\noindent
\textit{Proof of Theorem \ref{fund multi young}}.
For the sake of clarity, we divide the proof into three steps.\\
\textbf{Step 1}. Let $\m(\Om\times\square\times\R^m)$ be the set of finite real valued Radon measures 
on $\Om\times\square\times\R^m$ and let $\widehat{\nu}_h$ ($h=1,\ldots$) 
be the measure canonically associated with $w_h$,
\textit{i.e.},
\begin{equation*}
\widehat{\nu}_h(S)=
\int_{\Om}\chi_S\bigl(x, w_h(x)\bigr)\,dx \quad\text{for all} \;S\in\b(\Om\times\square\times\R^m).
\end{equation*}
Since $\widehat{\nu}_h$ is a bounded sequence of Radon measures, it has a weakly* convergent subsequence 
$\widehat{\nu}_{h_i}$. Let $\widehat{\nu}$ be the limit of $\widehat{\nu}_{h_i}$.
We claim that 
\begin{equation}\label{proiezione}
\widehat{\nu}\,(S\times\R^m)=\snorm{S} \quad\text{for all} \;S\in\b(\Om\times\square).
\end{equation}
We use a classical result about the convergence of \;$\widehat{\nu}_h(S)$\, when 
$S\subseteq\Om\times\square\times\R^m$ is an open or compact set (see \cite[Proposition~1.62]{Amb00}).

For $U\subseteq\Om$ and $V\subseteq\square$ open
\begin{equation}\label{pro uno}\begin{split}
\widehat{\nu}\,(U\times V\times\R^m)
\leq&\liminf_{i\rightarrow+\infty}\widehat{\nu}_{h_i}\,(U\times V\times\R^m)\\
=&\lim_{i\rightarrow+\infty}\int_{\Om}\chi_U(x)\chi_V\bigl(\langle\frac{x}{\e_{h_i}}\rangle\bigr)dx
=\snorm{U\times V},
\end{split}\end{equation}
by using Riemann-Lebesgue lemma (see \cite[Theorem 1.5]{Dac89}) in the last equality.

Denoting with $B_k$ the open ball with centre in the origin and radius $k\in\NP$, 
we have for $X\subseteq\Om$ and $Y\subseteq\square$ compact 
\begin{equation}\label{pro due}\begin{split}
\widehat{\nu}&\,(X\times Y\times\overline{B}_k)
\geq\limsup_{i\rightarrow+\infty}\widehat{\nu}_{h_i}\,(X\times Y\times\overline{B}_k)\\
&=\limsup_{i\rightarrow+\infty}\Bigl\vert\Bigl\lbrace x\in X \;:\; 
\langle\frac{x}{\e_{h_i}}\rangle\in Y \;\text{and}\; \snorm{u_{h_i}(x)}\leq k\Bigr\rbrace\Bigr\vert\\
&\geq\lim_{i\rightarrow+\infty}\int_{\Om}\chi_X(x)\chi_Y\bigl(\langle\frac{x}{\e_{h_i}}\rangle\bigr)dx
-\frac{\sup_i\norm{u_{h_i}}_{L^1}}{k}
=\snorm{X\times Y}-\frac{\sup_i\norm{u_{h_i}}_{L^1}}{k}.
\end{split}\end{equation}
From \eqref{pro uno} and \eqref{pro due}, by using inner and outer 
regularity of the measure $\widehat{\nu}$, we obtain equality~\eqref{proiezione}.

As a consequence, from disintegration theorem (see \cite[Theorem 1.5.1]{Ev90} and \cite[Theorem 2.28]{Amb00}), 
we can infer that there exists a map $\nu\in\Y(\Om\times\square,\R^m)$ such that 
\begin{equation*}
\int_{\Om\times\square\times\R^m}\phi\,d\,\widehat{\nu}
=\int_{\Om\times\square}\biggl(\int_{\R^m}\phi(x,y,\lambda)d\nu_{(x,y)}(\lambda)\biggr)dx\,dy 
\end{equation*}
for each continuous and bounded $\phi:\Om\times\square\times\R^m\rightarrow\R$.

\vspace{5pt}
\noindent
\textbf{Step 2}. We prove property (i). 
Note that by Remark \ref{misurab}, both sides of inequality \eqref{semicont} are well defined.
We first assume that there exists $b>0$ such that
\begin{equation}\label{inizio}
\W(x, y, \lambda)\leq b \quad\text{for all} \;(x, y,\lambda)\in\Om\times\square\times\R^m.
\end{equation}
By the admissibility condition, for every $\delta>0$ there exist a compact set $X_\delta\subseteq\Om$ and  
a compact set $Y_\delta\subseteq\square$ such that $\snorm{\Om\backslash X_\delta}\leq\delta$, 
$\snorm{\square\backslash Y_\delta}\leq\delta$ and $W\vert _{X_\delta\times Y_\delta\times\R^m}$ is continuous.
Let $\phi\in C(\Om\times\square\times\R^m)$ be an extension of $\W\vert_{X_\delta\times Y_\delta\times\R^m}$ 
such that $0\leq\phi(x, y, \lambda)\leq b$ for every $(x, y, \lambda)\in\Om\times\square\times\R^m$.
Since $\widehat{\nu}_{h_i}\rightharpoonup\widehat{\nu}$ weakly*,
defining \,$\overline{\phi}(x,y):=\int_{\R^m}\phi(x,y,\lambda)\;d\nu_{(x,y)}(\lambda)$, 
we have by the previous step
\begin{equation*}\begin{split}
\liminf_{i\rightarrow+\infty}\int_{\Om}\phi\bigl(x,w_{h_i}(x)\bigr)\,dx
=&\liminf_{i\rightarrow+\infty}\int_{\Om\times\square\times\R^m}\phi\,d\,\widehat{\nu}_{h_i}\\
\geq&\int_{\Om\times\square\times\R^m}\phi\,d\,\widehat{\nu}=
\int_{\Om\times\square}\overline{\phi}(x, y)\,dx\,dy.
\end{split}\end{equation*}

\vspace{8pt}
\noindent
We can write
\begin{equation*}\begin{split}
\int_{\Om\times\square}&\overline{\W}(x, y)\,dx\,dy -b\,\delta(1+\snorm{\Om})
\leq\int_{X_\delta\times Y_\delta}\overline{\W}(x, y)\,dx\,dy\\
&=\int_{X_\delta\times Y_\delta}\overline{\phi}(x, y)\,dx\,dy
\leq\int_{\Om\times\square}\overline{\phi}(x, y)\,dx\,dy
\leq\liminf_{i\rightarrow+\infty}\int_{\Om}\phi\bigl(x,w_{h_i}(x)\bigr)dx\\
&\leq\liminf_{i\rightarrow+\infty}\int_{\Om}
\Bigl[\chi_{X_\delta}(x)\chi_{Y_\delta}\bigl(\langle\frac{x}{\e_{h_i}}\rangle\bigr)
+\chi_{\Om\backslash X_\delta}(x)
+\chi_{\square\backslash Y_\delta}\bigl(\langle\frac{x}{\e_{h_i}}\rangle\bigr)\Bigl]
\phi\bigl(x,w_{h_i}(x)\bigr)dx\\
&\leq\liminf_{i\rightarrow+\infty}\int_{X_\delta}
\chi_{Y_\delta}\bigl(\langle\frac{x}{\e_{h_i}}\rangle\bigr)\W\bigl(x,w_{h_i}(x)\bigr)dx+b\,\delta
+b\lim_{i\rightarrow+\infty}\int_{\Om}\chi_{\square\backslash Y_\delta}\bigl(\langle\frac{x}{\e_{h_i}}\rangle\bigr)dx\\
&\leq\liminf_{i\rightarrow+\infty}\int_{\Om}\W\bigl(x,w_{h_i}(x)\bigr)dx
+b\,\delta(1+\snorm{\Om}).
\end{split}\end{equation*}

\vspace{2pt}
\noindent
Being $\delta>0$ arbitrary, inequality \eqref{semicont} follows.

In order to remove assumptions \eqref{inizio} we consider, for $k\in\NP$, the functions 
\begin{equation*}
\W_k(x,y,\lambda):=\mathrm{min}\bigl\lbrace k, \W(x,y,\lambda)\bigr\rbrace.
\end{equation*}

\vspace{3pt}
\noindent
By applying the first part of the step, we have
\begin{equation*}
\liminf_{i\rightarrow+\infty}\int_{\Om}\W\bigl(x,w_{h_i}(x)\bigr)dx\geq
\liminf_{i\rightarrow+\infty}\int_{\Om}\W_k\bigl(x,w_{h_i}(x)\bigr)dx\geq
\int_{\Om\times\square}\overline{\W}_k(x, y)\,dx\,dy.
\end{equation*}

\vspace{3pt}
\noindent
By noting that $\W_k$ is increasing and that $\W_k(x, y, \cdot)\rightarrow\W(x, y, \cdot)$ 
\textit{a.e.} in $\R^m$ for every fixed $(x,y)\in\Om\times\square$, we deduce from the monotone convergence 
theorem that $\overline{\W}_k\rightarrow\overline{\W}$ \textit{a.e.} in $\Om\times\square$. 
The sequence $\overline{\W}_k$ is increasing so, again from monotone convergence theorem, 
\begin{equation*}
\int_{\Om\times\square}\overline{\W}_k(x, y)\,dx\,dy\xrightarrow{k\rightarrow\infty} 
\int_{\Om\times\square}\overline{\W}(x, y)\,dx\,dy.
\end{equation*}

\vspace{5pt}
\noindent
\textbf{Step 3}.
We prove property (ii). Let $\W^+:=\mathrm{max}\lbrace 0, \W\rbrace$ and $\W^-:=\mathrm{max}\lbrace 0, -\W\rbrace$,
so that $\W=\W^+-\W^-$. Both the sequences $\W^+(\cdot,w_{h}(\cdot))$ and $\W^-(\cdot,w_{h}(\cdot))$ are
equi-integrable, so it is enough to prove equality \eqref{cont} when $\W$ is non-negative.

If $\W$ is bounded from above by a constant $b$, then \eqref{cont} follows by applying \eqref{semicont} 
to $\W$ and $b-\W$.
For general non-negative $\W$, by an equivalent characterization of equi-integrability, 
for each $\eta>0$ there exists $k\in\NP$ such that 
$\mathrm{sup}_h\int_{\lbrace x\,:\,\W(x,w_h(x))\geq k\rbrace}\W\bigl(x,w_h(x)\bigr)\,dx<\eta$. Hence
\begin{equation*}\begin{split}
\limsup_{i\rightarrow+\infty}&\int_{\Om}\W\bigl(x,w_{h_i}(x)\bigr)dx-\eta
\leq\lim_{i\rightarrow+\infty}\int_{\Om}\W_k\bigl(x,w_{h_i}(x)\bigr)dx\\
=&\int_{\Om\times\square}\overline{\W}_k(x,y)\,dx\,dy 
\leq\int_{\Om\times\square}\overline{\W}(x,y)\,dx\,dy.
\end{split}\end{equation*}

\vspace{2pt}
\noindent
Being $\eta>0$ arbitrary, the previous inequality completes the proof.
\vspace{-10pt}
\begin{flushright}\qedsymbol\end{flushright}

\vspace{6pt}
\begin{prop}\label{local supporti}
Let $u_h$ be a bounded sequence in $L^1(\Om,\R^m)$ generating a two-scale Young measure $\nu$.
If $\lbrace\mathfrak{A}_k\rbrace_{k=1,\ldots,n}$ is a family of closed subsets of $\R^m$ and 
$\lbrace P_k\rbrace_{k=1,\ldots,n}$ is a measurable partition of the unit cell $\square$, then 
the following conditions are equivalent:

\vspace{5pt}
\begin{itemize}
\item[i)] for all \;$k\in\lbrace1,\ldots,n\rbrace$ and for \textit{a.e.} $(x,y)\in \Om\times P_k$ 
\;$\mathrm{supp}\,\nu_{(x,y)}\subseteq\mathfrak{A}_k$;

\vspace{5pt}
\item[ii)] $\D\bigl(\langle\frac{x}{\e_h}\rangle,u_h(x)\bigr)\rightarrow 0$ in measure, where 
$\D(y,\lambda):=\sum_{k=1}^{n}\chi_{P_k}(y)\mathrm{dist}(\lambda,\mathfrak{A}_k)$.
\end{itemize}
\end{prop}

\vspace{2pt}
\begin{proof}
(i)$\Rightarrow$(ii) 
Fixed $\eta>0$, we define $\W(y,\lambda):=\mathrm{min}\lbrace\eta,\D(y,\lambda)\rbrace$.
By Remark \ref{scorza}, $\W$ is an admissible integrand and therefore, 
by Theorem \ref{fund multi young},
\begin{equation*}
\lim_{h\rightarrow+\infty}\int_{\Om}\W\Bigl(\langle\frac{x}{\e_h}\rangle, u_h(x)\Bigr)dx
=\int_{\Om\times\square}\biggl(\int_{\R^m}\W(y,\lambda)\;d\nu_{(x,y)}(\lambda)\biggr)dx\,dy=0.
\end{equation*}
Since $\eta>0$ is arbitrary, this proves that 
$\D\bigl(\langle\frac{x}{\e_h}\rangle,u_h(x)\bigr)\rightarrow 0$ in measure.

(ii)$\Rightarrow$(i) 
Let $\varphi\in C_k:=\lbrace\varphi\in C_c(\R^m) \,:\, \mathrm{supp}\,\varphi\subseteq\R^m\setminus\mathfrak{A}_k\rbrace$
and $\eta>0$ such that $\varphi(\lambda)=0$ if $\mathrm{dist}(\lambda,\mathfrak{A}_k)\leq\eta$.
Taking 
\begin{equation*}
S_h:=\Bigl\lbrace x\in\Om \,:\, 
\langle\frac{x}{\e_h}\rangle\in P_k \;\text{and}\; \mathrm{dist}\bigl(u_h(x),\mathfrak{A}_k\bigr)>\eta\Bigr\rbrace, 
\end{equation*}
it follows from the hypothesis that $\snorm{S_h}\rightarrow0$. 
Given $X\subseteq\Om$ and $Y\subseteq P_k$ measurable, we define $\W(x,y,\lambda):=\chi_X(x)\chi_Y(y)\varphi(\lambda)$. 
By Lusin theorem, $\W$ is an admissible integrand and therefore, 
by Theorem \ref{fund multi young},
\begin{equation*}
\int_{X\times Y}\biggl(\int_{\R^m}\varphi(\lambda)\;d\nu_{(x,y)}(\lambda)\biggr)dx\,dy
=\lim_{h\rightarrow+\infty}\int_{S_h}\W\Bigl(x,\langle\frac{x}{\e_h}\rangle, u_h(x)\Bigr)dx=0.
\end{equation*}
Since $X$ and $Y$ are arbitrary and $C_k$ is separable, it remains proved that for \textit{a.e.} 
$(x,y)\in\Omega\times P_k$ we have $\int_{\R^m}\varphi\;d\nu_{(x,y)}=0$ for all 
$\varphi\in C_k$ or equivalently $\mathrm{supp}\,\nu_{(x,y)}\subseteq\mathfrak{A}_k$.
\\
\end{proof}

\vspace{2pt}
\begin{defi}
A map \,$\nu\in\Y(\Om\times\square,\M^{d\times s})$ is a \emph{two-scale gradient Young measure} if
there exists a sequence $u_h$ in $W^{1,1}(\Om, \R^s)$ such that $\g u_h$ is bounded in $L^1(\Om,\M^{d\times s})$ 
and generates $\nu$ in two-scale.
For a complete characterization we refer to \cite{BBS06} (see also \cite{Ped05}).
\end{defi}
%%%%%%%%%%%%%%%%%%%%%%%%%%%%%%%%%%%%%%%%%%%%%%%%%%%%%%%%%%%%%%%%%%%%%%%%%%%%%%%%%%%%%%%%%%%%%%%%%%%
%%%%%%%%%%%%%%%%%%%%%%%%%%%%%%%%%%%%%%%%%%%%%%%%%%%%%%%%%%%%%%%%%%%%%%%%%%%%%%%%%%%%%%%%%%%%%%%%%%%
\section{Gamma-convergence}
\noindent
In this section we give a proof of Theorem \ref{zero energy}.
First we recall the definition of $\Gamma$-convergence,
referring to \cite{Br98} and \cite{Dal93} for a comprehensive treatment. 

\begin{defi}
Let $(U,\tau)$ be a topological space satisfying the first countability axiom and
$E_h$, $E$ functionals from $U$ to $[-\infty, +\infty]$;
we say that $E$ is the $\Gamma(\tau)$-limit of the sequence $E_h$, or that $E_h$ $\Gamma(\tau)$-converges to $E$,
and write
\begin{equation*}
E=\Gamma(\tau)\hbox{-}\negthinspace\lim_{h\rightarrow+\infty}E_h,
\end{equation*}
if for every $u\in U$ the following conditions are satisfied:

\begin{equation*}
\hspace*{-5pt}E(u)\leq\inf\Bigl\lbrace\liminf_{h\rightarrow+\infty}E_h(u_h) : u_h\xrightarrow{\tau}u\Bigr\rbrace
\end{equation*}
and
\begin{equation*}
E(u)\geq\inf\Bigl\lbrace\limsup_{h\rightarrow+\infty}E_h(u_h) : u_h\xrightarrow{\tau}u\Bigr\rbrace.
\end{equation*}
\end{defi}

\vspace{8pt}
We can extend the definition of $\Gamma$-convergence to families depending on a parameter $\e>0$.
\begin{defi}
For every $\e>0$, let ${E_\e}$ be a functional from $U$ to $[-\infty, +\infty]$.
We say that $E$ is the $\Gamma(\tau)$-limit of the family $E_\e$ as $\e\rightarrow0^+$, and write
\begin{equation*}
E=\Gamma(\tau)\hbox{-}\negthinspace\lim_{\e\rightarrow0^+}E_\e,
\end{equation*}
if we have for every sequence $\e_h\rightarrow0^+$
\begin{equation*}
E=\Gamma(\tau)\hbox{-}\negthinspace\lim_{h\rightarrow+\infty}E_{\e_h}.
\end{equation*}
\end{defi}

\vspace{10pt}
We are interested in $\Gamma$-convergence for periodic homogenization of integral functionals.
In the following, the space $L^p(\Om,\R^s)$, $p\in(1,+\infty)$, is endowed with the strong topology. 
We consider the family of functionals $E_\e:L^p(\Om,\R^s)\rightarrow[0, +\infty]$ defined by
\begin{equation*}
E_\e (u):=
\begin{cases}
\displaystyle\int_\Om \W\Bigl(\langle\frac{x}{\e}\rangle,\g u(x)\Bigr)dx
& \text{if} \quad u\in W^{1,p}(\Om,\R^s),\\
&\\
+\infty & \text{otherwise,}
\end{cases}
\end{equation*}
where $\W:\square\times\M^{s\times d}\rightarrow[0, +\infty)$ is a Carath\'eodory
function satisfying coerciveness and growth conditions of order $p$ as in \eqref{crescita}.

\vspace{4pt}
The next theorem is a standard result (see \cite[Theorem 14.5]{Br98}).
\begin{teo}\label{Gamma}
The family $E_\e$ \,$\Gamma(L^p)$-converges and its $\Gamma(L^p)$-limit 
\,$E_{hom}:W^{1,p}(\Om,\R^s)\rightarrow[0, +\infty]$ can be written as 
\begin{equation*}
E_{hom} (u)=
\begin{cases}
\displaystyle\int_\Om \W_{hom}\bigl(\g u(x)\bigr)dx
& \text{if} \quad u\in W^{1,p}(\Om,\R^s),\\
&\\
+\infty & \text{otherwise,}
\end{cases}
\end{equation*}

\vspace{5pt}
\noindent
where $\W_{hom}:\M^{s\times d}\rightarrow[0, +\infty)$ is a suitable quasiconvex function 
depending on $\W$.
\end{teo}

\begin{defi}
We call $\W_{hom}$ the \emph{homogenized integrand} related to $\W$.
\end{defi}

\vspace{9pt}
In the proof of Theorem \ref{zero energy}, we will use the following lemma, 
which can be derived by \cite[Theorem~4]{Mul99} (see also \cite[Lemma~3.1]{Zha92}).

\begin{lem}\label{regularity}
Assume that $\mathfrak{A}\subseteq\M^{s\times d}$ is a compact set.
Let $A\in\M^{s\times d}$ and let $u_h$ be a sequence converging weakly in $W^{1,1}(\Om,\R^s)$ to $A\,x$
and such that $\mathrm{dist}(\g u_h,\mathfrak{A})\rightarrow0$ in measure. 
Then there exists a sequence $v_h$ in $W^{1,\infty}(\Om,\R^s)$ such that
\begin{itemize}
\item[i)] $v_h\rightharpoonup A\,x$ weakly* in $W^{1,\infty}(\Om,\R^s)$;
\vspace{4pt}
\item[ii)]  $\snorm{\lbrace \g v_h\neq\g u_h\rbrace}\rightarrow0$.
\end{itemize}
\end{lem}

\vspace{8pt}
\noindent
\textit{Proof of Theorem \ref{zero energy}}.
We begin to prove the inclusion $\mathfrak{A}_{hom}\subseteq\W_{hom}^{-1}(0)$.
Note that, by coerciveness hypothesis, the sets $\mathfrak{A}_1$, ..., $\mathfrak{A}_n$
are compact. Fix $A\in\mathfrak{A}_{hom}$ and $\e_h\rightarrow0^+$. 
By definition of $\mathfrak{A}_{hom}$, there exists a sequence $u_h$ such that 
$u_h\rightharpoonup A\,x$ weakly* in $W^{1,\infty}(\Om,\R^s)$ and
$\D\bigl(\langle\frac{x}{\e_h}\rangle,\g u_h(x)\bigr)\rightarrow 0$ in measure, 
where $\D(y,\Lambda):=\sum_{k=1}^{n}\chi_{P_k}(y)\mathrm{dist}(\Lambda,\mathfrak{A}_k)$.
For a suitable subsequence $h_i$, $\g u_{h_i}$ generates a two-scale Young measure
$\nu$ with respect to $\e_{h_i}$.

By Proposition \ref{local supporti}, for all $k\in\lbrace1,\ldots,n\rbrace$ and
for \textit{a.e.} $(x,y)\in \Om\times P_k$ \,we have \,$\mathrm{supp}\,\nu_{(x,y)}\subseteq\mathfrak{A}_k$, 
therefore $\int_{\M^{s\times d}}\W(y,\Lambda)\;d\nu_{(x,y)}(\Lambda)=0$.
Observed that $\W$ is an admissible integrand, we obtain by Theorem \ref{fund multi young}
\begin{equation*}
\lim_{i\rightarrow+\infty}\int_{\Om}\W\Bigl(\langle\frac{x}{\e_{h_i}}\rangle,\g u_{h_i}(x)\Bigr)dx
=\int_{\Om\times\square}\biggl(\int_{\M^{s\times d}}\W(y,\Lambda)\;d\nu_{(x,y)}(\Lambda)\biggr)dx\,dy=0.
\end{equation*}
By the definition of $\Gamma(L^p)$-limit, $\W_{hom}(A)=0$. 

\vspace{2pt}
We prove now the opposite inclusion.
Fix $A\in\W_{hom}^{-1}(0)$ and $\e_h\rightarrow0^+$. 
By definition of $\W_{hom}$, there exists a sequence $u_h$ such that $u_h\rightarrow A\,x$ 
strongly in $L^p(\Om,\R^s)$ and $E_{\e_h}(u_h)\rightarrow 0$. 
For $h$ large enough, $E_{\e_h}(u_h)$ is finite and therefore
$u_h\in W^{1,p}(\Om,\R^s)$. Thanks to the $p$-coerciveness hypothesis on $\W$, we can 
suppose that $u_h\rightharpoonup A\,x$ weakly in $W^{1,p}(\Om,\R^s)$.

We claim that $\D\bigl(\langle\frac{x}{\e_h}\rangle,\g u_h(x)\bigr)\rightarrow 0$ in measure.
If not, there exist $\delta,\eta>0$ and a subsequence~$u_{h_i}$ such that
\begin{equation}\label{assurdo}
\inf_i\,\Bigl\vert\Bigl\lbrace x\in\Om \;:\; 
\D\bigl(\langle\frac{x}{\e_{h_i}}\rangle,\g u_{h_i}(x)\bigr)>\eta\Bigr\rbrace\Bigr\vert>\delta.
\end{equation}
Refining the subsequence $h_i$ if necessary, we can suppose that $\g u_{h_i}$ generates a two-scale 
Young measure $\nu$ with respect to $\e_{h_i}$. By Theorem \ref{fund multi young} 
\begin{equation*}
0=\lim_{i\rightarrow+\infty}\int_{\Om}\W\Bigl(\langle\frac{x}{\e_{h_i}}\rangle,\g u_{h_i}(x)\Bigr)dx
\geq\int_{\Om\times\square}\biggl(\int_{\M^{s\times d}}\W(y,\Lambda)\;d\nu_{(x,y)}(\Lambda)\biggr)dx\,dy
\end{equation*}

\vspace{3pt}
\noindent
so that for all $k\in\lbrace1,\ldots,n\rbrace$ and for \textit{a.e.} $(x,y)\in\Om\times P_k$
\,we have \,$\int_{\M^{s\times d}}\W_k(\Lambda)\;d\nu_{(x,y)}(\Lambda)=0$. 
Since $\W_k$ is continuous and non-negative,
$\mathrm{supp}\,\nu_{(x,y)}\subseteq\mathfrak{A}_k$ for \textit{a.e.} $(x,y)\in\Om\times P_k$.
Therefore by Proposition~\ref{local supporti} 
\;$\D\bigl(\langle\frac{x}{\e_{h_i}}\rangle,\g u_{h_i}(x)\bigr)\rightarrow 0$ in measure,
in contradiction with \eqref{assurdo}.

Finally, noticed that $\mathrm{dist}\bigl(\g u_h(x),\bigcup_{k=1}^n\mathfrak{A}_k\bigr)\rightarrow 0$ in measure,
we can apply Lemma \ref{regularity} to infer the existence of a sequence $v_h$ such that $v_h\rightharpoonup A\,x$ 
weakly* in $W^{1,\infty}(\Om,\R^s)$ and $\snorm{\lbrace\g v_h\neg\neq\neg\g u_h\rbrace}\rightarrow0$.
In particular, $\D\bigl(\langle\frac{x}{\e_h}\rangle,\g v_h(x)\bigr)\rightarrow0$ in measure.
\vspace{-6pt}
\begin{flushright}\qedsymbol\end{flushright}
%%%%%%%%%%%%%%%%%%%%%%%%%%%%%%%%%%%%%%%%%%%%%%%%%%%%%%%%%%%%%%%%%%%%%%%%%%%%%%%%%%%%%%%%%%%%%%%%%%%
%%%%%%%%%%%%%%%%%%%%%%%%%%%%%%%%%%%%%%%%%%%%%%%%%%%%%%%%%%%%%%%%%%%%%%%%%%%%%%%%%%%%%%%%%%%%%%%%%%%
\section{Polyconvexity}
\noindent
In this section we recall some of the definitions and the results related to polyconvexity. 
General references are \cite{Dac89}, \cite{Dac99} and \cite{Mul99bis}.
In the following we always assume that $d,s\geq2$ since otherwise polyconvexity agrees with 
ordinary convexity.

\begin{defi}
A function $\W:\M^{s\times d}\rightarrow\R$ is said to be
\emph{polyconvex} if there is a convex function $\V:\R^{\tau(d,s)}\rightarrow\R$ such that
$\W(\Lambda)=\V(\widehat{\Lambda})$ for all $\Lambda\in\M^{s\times d}$. Here $\widehat{\Lambda}$ 
denotes the list of all minors (subdeterminants) of $\Lambda$ and $\tau(d,s)=\binom{d+s}{d}-1$. 
We can identify $\widehat{\Lambda}$ with a point of $\R^{\tau(d,s)}$. 
In the simple case $d=s=2$\; we have $\widehat{\Lambda}=(\Lambda,\mathrm{det}(\Lambda))\in\R^5$. 
\end{defi}
\begin{rem}
The maximum of two polyconvex functions is still polyconvex.
\end{rem}

\vspace{4pt}
\begin{defi} 
We say that a set $\mathfrak{A}\subseteq\M^{s\times d}$ is \emph{polyconvex} if there exists a convex set 
$\mathfrak{C}\subseteq\R^{\tau(d,s)}$ such that
\begin{equation}\label{polyconv}
\mathfrak{A}=\Bigl\lbrace A\in\M^{s\times d} \,:\, \widehat{A}\in\mathfrak{C}\Bigr\rbrace.
\end{equation}
\end{defi}
\begin{defi}
The \emph{polyconvex hull} of a set $\mathfrak{A}\subseteq\M^{s\times d}$ is the smallest polyconvex
set containing $\mathfrak{A}$ and is denoted by $\mathfrak{A}^{pc}$.
\end{defi}

\vspace{6pt}
\begin{lem}\label{esistenza azzeranti}
Let $\mathfrak{A}$ be a compact set of $\M^{s\times d}$ and let \;$p\geq\mathrm{max}\lbrace d,s\rbrace$. 
If $\mathfrak{A}$ is polyconvex, then there exists a polyconvex function $\W:\M^{s\times d}\rightarrow[0,+\infty)$ satisfying coerciveness and growth conditions of order $p$ such that 
\begin{equation*}
\mathfrak{A}=\W^{-1}(0).
\end{equation*}
\end{lem}
\begin{proof}
Let $\mathfrak{C}\subseteq\R^{\tau(d,s)}$ be a compact and convex set such that \eqref{polyconv} holds.
Consider the functions
\begin{equation*}
\W_1(\Lambda):=\mathrm{dist}^p(\Lambda,\mathfrak{A}^{co}) \quad\text{and}\quad
\W_2(\Lambda):=\mathrm{dist}^{\frac{p}{q}}(\widehat{\Lambda},\mathfrak{C}),
\end{equation*}
where $q=\mathrm{max}\lbrace d,s\rbrace$ and $co$ denotes the convex hull.
Both are polyconvex and with $p$-growth, moreover $\W_1$ is $p$-coercive and $\W_2^{-1}(0)=\mathfrak{A}$.
The function $\W:=\mathrm{max}\lbrace\W_1,\W_2\rbrace$ does the job. 
\\
\end{proof}

%%%%%%%%%%%%%%%%%%%%%%%%%%%%%%%%%%%%%%%%%%%%%%%%%%%%%%%%%%%%%%%%%%%%%%%%%%%%%%%%%%%%%%%%%%%%%%%%%%%
%%%%%%%%%%%%%%%%%%%%%%%%%%%%%%%%%%%%%%%%%%%%%%%%%%%%%%%%%%%%%%%%%%%%%%%%%%%%%%%%%%%%%%%%%%%%%%%%%%%
\section{Proof of Theorem \ref{main}}
\noindent
This section is devoted to the proof of Theorem \ref{main}.
Assuming in the following that $d=2$, we begin by a technical lemma. 
\begin{lem}\label{sverak split}
Let $\mathfrak{A}_{hom}$ be the homogenized set related to $\lbrace(\mathfrak{A}_k,P_k)\rbrace_{k=1,2}$, 
where $\mathfrak{A}_1$, $\mathfrak{A}_2$, $P_1$ and $P_2$ are defined as in Theorem \ref{main}.
If $A\in\mathfrak{A}_{hom}$, then $A-\dig\left(\frac{1}{2},1\right)$ is not positive definite.
\end{lem}

This lemma will be proved by using a remarkable result, due to \v Sver\'ak.
\begin{lem}\label{sverak}\textnormal{(See \cite[Corollary 1]{Sver92} and also \cite{FarX03})}.
Let $\V:\M^{2\times2}\nek\rightarrow\nek[0,+\infty)$ be the continuous function defined by
\begin{equation}\label{funz sverak}
\V(\Lambda)=
\begin{cases}
\displaystyle\mathrm{det}\,\mathcal{E}(\Lambda)
& \text{if} \;\Lambda \;\text{is positive definite},\\
0 & \textrm{otherwise,}
\end{cases}
\end{equation}
where $\mathcal{E}(\Lambda)$ denotes the symmetric part of $\Lambda$.
Let $v_h\in W^{1,\infty}(\Om,\R^2)$ be a sequence such that
\begin{itemize}
\item[i)]   $v_h\rightharpoonup A\,x$ weakly* in $W^{1,\infty}(\Om,\R^2)$;
\vspace{4pt}
\item[ii)]  $\mathrm{dist}(\g v_h,\M^{2\times2}_{sym})\rightarrow0$ in measure;
\vspace{4pt}
\item[ii)]  $\g v_h$ generates a Young measure $\mu\in\Y(\Om,\M^{2\times2})$.
\end{itemize}
\vspace{4pt}
Then $\int_{\M^{2\times2}}\V(\Lambda)\,d\mu_x(\Lambda)\geq\V(A)$ for \textit{a.e.} $x\in\Om$.
\end{lem}

\vspace{7pt}
\noindent
\textit{Proof of Lemma \ref{sverak split}}.
Define as usual 
$\D(y,\Lambda)
:=\chi_{P_1}(y)\mathrm{dist}(\Lambda,\mathfrak{A}_1)+\chi_{P_2}(y)\mathrm{dist}(\Lambda,\mathfrak{A}_2)$.
Let $\e_h\nec\rightarrow0^+$ and let $u_h$ be a sequence such that $u_h\rightharpoonup A\,x$ weakly* in
$W^{1,\infty}(\Om,\R^2)$ and $\D\bigl(\langle\frac{x}{\e_h}\rangle,\g u_h(x)\bigr)\rightarrow 0$ in measure.
Passing to a subsequence if necessary, we can assume that $\g u_h$ generates a two-scale Young measure $\nu$.
By Proposition \ref{local supporti}, we have \,$\mathrm{supp}\,\nu_{(x,y)}\subseteq\mathfrak{A}_k$ 
for \textit{a.e.} $(x,y)\in \Om\times P_k$ \;($k=1,2$).

Let observe now that $\dig\left(\frac{1}{2},1\right)=\frac{1}{2}\left[\dig{(-2,1)}+\dig{(3,1)}\right]$ and that
the matrices belonging to $\mathfrak{A}_1-\dig{(-2,1)}$ and $\mathfrak{A}_2-\dig{(3,1)}$ are not positive definite.
In particular, if $\V$ is defined as in~\eqref{funz sverak}, then
\begin{equation}\label{spezzamento}\begin{split}
&\int_{\M^{2\times2}}\V\bigl(\Lambda-\dig{(-2,1)}\bigr)d\nu_{(x,y)}(\Lambda)=0 \quad\text{for \textit{a.e.}}(x,y)\in\Om\times P_1\\
&\int_{\M^{2\times2}}\V\bigl(\Lambda-\dig{(3,1)}\bigr)d\nu_{(x,y)}(\Lambda)=0 \;\;\;\quad
                                                                       \text{for \textit{a.e.}}(x,y)\in\Om\times P_2.
\end{split}\end{equation}

Let $v_h\in W^{1,\infty}(\square,\R^2)$ be the sequence defined by
\begin{equation}\label{saltante}\begin{split}
v_h(x)_1&:=u_h(x)_1-\frac{5}{2}\,\e_h\Bigl\vert\langle\frac{x_1}{\e_h}\rangle-\frac{1}{2}\Bigr\vert
        -\frac{1}{2}\,x_1\\
v_h(x)_2&:=u_h(x)_2-x_2.
\end{split}\end{equation}

\vspace{3pt}
\noindent
It is trivial to check that
$\g v_h(x)=\g u_h(x)-\chi_{P_1}(\langle\frac{x}{\e_h}\rangle)\dig(-2,1)
                    -\chi_{P_2}(\langle\frac{x}{\e_h}\rangle)\dig(3,1)$
and $\mathrm{dist}\bigl(\g v_h(x),\M^{2\times2}_{sym}\bigr)\rightarrow0$ in measure.
Moreover by Riemann-Lebesgue lemma follows that
$v_h\rightharpoonup \bigl(A-\dig(\frac{1}{2},1)\bigr)x$ weakly* in $W^{1,\infty}(\Om,\R^2)$.
Passing to a subsequence if necessary, we can suppose that $\g v_h$ generates a Young measure $\mu$.

Consider the function
\begin{equation*}
\W(y,\Lambda)
:=\chi_{P_1}(y)\,\V\bigl(\Lambda-\dig{(-2,1)}\bigr)+\chi_{P_2}(y)\,\V\bigl(\Lambda-\dig{(3,1)}\bigr).
\end{equation*}
By construction of the sequence $v_h$, we have
$\W\Bigl(\langle\frac{x}{\e_h}\rangle,\g u_h(x)\Bigr)=\V\bigl(\g v_h(x)\bigr)$.
Moreover, by Theorem \ref{fund multi young} and equalities \eqref{spezzamento}
\begin{equation*}
\lim_{h\rightarrow+\infty}\int_{\Om}\W\Bigl(\langle\frac{x}{\e_h}\rangle,\g u_h(x)\Bigr)dx
=\int_{\Om\times\square}\biggl(\int_{\M^{2\times2}}\W(y,\Lambda)\;d\nu_{(x,y)}(\Lambda)\biggr)dx\,dy=0\,
\end{equation*}
while by Theorem \ref{fund young} and Lemma \ref{sverak}
\begin{equation*}
\lim_{h\rightarrow+\infty}\int_{\Om}\V\bigl(\g v_h(x)\bigr)dx
=\int_{\Om}\biggl(\int_{\M^{2\times2}}\V(\Lambda)\;d\mu_x(\Lambda)\biggr)dx
\geq\int_{\Om}\V\Bigl(A-\dig\Bigl(\frac{1}{2},1\Bigl)\Bigr)dx.
\end{equation*}

\vspace{4pt}
\noindent
Since the last integrand vanishes, $\V\bigl(A-\dig\bigl(\frac{1}{2},1\bigl)\bigr)=0$, 
\textit{i.e.}, the matrix $A-\dig\left(\frac{1}{2},1\right)$ is not positive definite.
\vspace{-6pt}
\begin{flushright}\qedsymbol\end{flushright}

\begin{figure}
\label{disegno}
\begin{center}
\psfrag{0}{$O$}
\psfrag{1}{$A^1_1$}
\psfrag{2}{$A^2_1$}
\psfrag{3}{$B^1$}
\psfrag{4}{$B^2$}
\psfrag{5}{$A^1_2$}
\psfrag{6}{$A^2_2$}
\psfrag{7}{$(-2,1)$}
\psfrag{8}{$(\frac{1}{2},1)$}
\psfrag{9}{$(3,1)$}
\includegraphics[width=12cm]{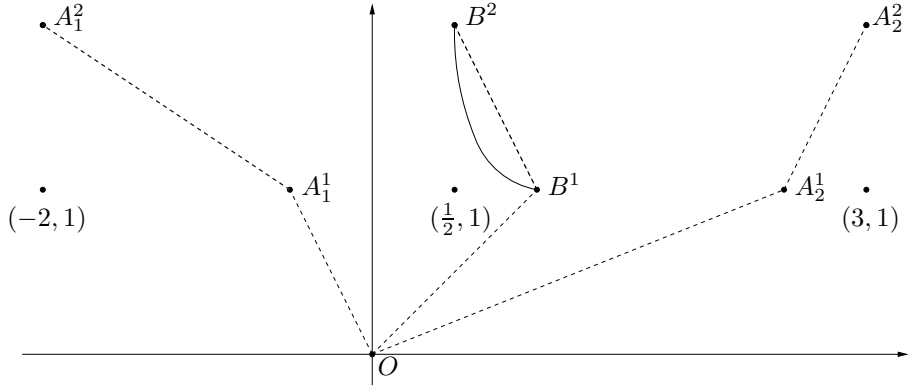}
\caption{A representation of $\mathfrak{A}_1$, $\mathfrak{A}_2$ and $\mathfrak{B}$ in $\R^2$, identified with 
the set of the diagonal matrices. The arc joining $B^1$ and $B^2$ represents the set
$\lbrace\dig(b_1(t),b_2(t)) \,:\, t\in(0,1)\rbrace$.}
\end{center}
\end{figure}

\vspace{6pt}
\noindent
\textit{Proof of Theorem \ref{main}}.
Let $\mathfrak{C}=\lbrace C^1,C^2,C^3\rbrace\subseteq\M^{2\times2}$ and let $\psi$ be the symmetric bilinear form on 
$\M^{2\times2}$ determined by $\mathrm{det}\,\Lambda=\frac{1}{2}\psi(\Lambda,\Lambda)$. 
We can write for all $(t_1,t_2,t_3)\in\R^3$
\begin{equation*}\begin{split}
\sum_{k,j=1}^3 t_kt_j\mathrm{det}(C_k-C_j)
&=\sum_{k,j=1}^3 t_kt_j\psi(C_k,C_k)-\sum_{k,j=1}^3 t_kt_j\psi(C_k,C_j)\\
&=2\sum_{k,j=1}^3 t_kt_j\mathrm{det}(C_k)-2\mathrm{det}\Bigl(\sum_{k=1}^3t_kC_k\Bigr).
\end{split}\end{equation*}
In particular
\begin{equation*}\begin{split}
\mathfrak{C}^{pc}&=\biggl\lbrace\sum_{k=1}^3t_kC_k \;:\; t_k\geq0, \;\sum_{k=1}^3t_k=1 \;\;\text{and}
                  \;\;\mathrm{det}\Bigl(\sum_{k=1}^3t_kC_k\Bigr)=\sum_{k=1}^3t_k\,\mathrm{det}\,C_k \biggr\rbrace\\
&=\biggl\lbrace\sum_{k=1}^3t_kC_k \;:\; t_k\geq0, \;\sum_{k=1}^3t_k=1 \;\;\text{and}
                  \;\;\sum_{k,j=1}^3t_kt_j\,\mathrm{det}(C_k-C_j)=0 \biggr\rbrace.
\end{split}\end{equation*}

\vspace{5pt}
\noindent
By the formula above, since $\mathrm{det}\,(\Lambda-\Lambda')$ does not change sign on $\mathfrak{A}_1\times\mathfrak{A}_1$
and $\mathfrak{A}_2\times\mathfrak{A}_2$, it follows that
\begin{equation*}
\mathfrak{A}_1=\mathfrak{A}_1^{pc} \quad\text{and}\quad \mathfrak{A}_2=\mathfrak{A}_2^{pc}.
\end{equation*}

\vspace{3pt}
\noindent
Moreover, an elementary computation shows that
$\mathfrak{B}^{pc}=\mathfrak{B}\cup\lbrace\dig(b_1(t),b_2(t)) \,:\, t\in(0,1)\rbrace$,
where $b_1$, $b_2$ are defined by
\begin{equation*}
b_1(t):=\frac{-3t+2+\sqrt{9t^2-4t+4}}{4} \quad \text{and} \quad b_2(t):=\frac{3t+2+\sqrt{9t^2-4t+4}}{4}.
\end{equation*}

\vspace{5pt}
By Lemma \ref{esistenza azzeranti}, there exist two polyconvex functions \,$\W_1,\W_2:\M^{2\times2}\nek\rightarrow\nek[0,+\infty)$ $p$-coercive and with $p$-growth such that 
\begin{equation*}
\mathfrak{A}_1=\W_1^{-1}(0) \quad\text{and}\quad \mathfrak{A}_2=\W_2^{-1}(0).
\end{equation*}

Fixed $\e_h\nec\rightarrow0^+$ and $k\in\lbrace1,2\rbrace$, it is easy to build as in \eqref{saltante} a sequence 
such that $u_h\rightharpoonup B^kx$ weakly* in $W^{1,\infty}(\Om,\R^2)$ and
$\g u_h(x)=\chi_{P_1}(\langle\frac{x}{\e_h}\rangle)A_1^k+\chi_{P_2}(\langle\frac{x}{\e_h}\rangle)A_2^k$.
Therefore
\begin{equation*}
\mathfrak{B}\subseteq\mathfrak{A}_{hom}.
\end{equation*}

\noindent
On the other hand, noted that $b_1(t)>\frac{1}{2}$ and $b_2(t)>1$ for all $t\in(0,1)$,  
we have from Lemma~\ref{sverak split}
\begin{equation*}
\mathfrak{B}^{pc}\nsubseteq\mathfrak{A}_{hom}.
\end{equation*} 

Finally, the homogenized integrand $\W_{hom}$ related to $\W(y,\Lambda):=\chi_{P_1}(y)\,\W_1(\Lambda)+\chi_{P_2}(y)\,\W_2(\Lambda)$
cannot be polyconvex because by Theorem \ref{zero energy} \;$\mathfrak{A}_{hom}=\W_{hom}^{-1}(0)$.
\vspace{-6pt}
\begin{flushright}\qedsymbol\end{flushright}

\vspace{6pt}
\begin{rem}
The polyconvex hull of a compact set $\mathfrak{C}\subseteq\M^{2\times2}$ can be characterized as
\begin{equation*}
\mathfrak{C}^{pc}=\biggl\lbrace\int_{\M^{2\times2}}\Lambda\,d\nu \;:\; \nu\in\p(\M^{2\times2})\,, 
                         \;\;\mathrm{supp}\,\nu\subseteq\mathfrak{C} \;\;\text{and}\; 
                     \int_{\M^{2\times2}}\mathrm{det}\,\Lambda\,d\nu=\mathrm{det}\int_{\M^{2\times2}}\Lambda\,d\nu \biggr\rbrace.
\end{equation*}
Thanks to this characterization, the polyconvexity of $\mathfrak{A}_1$ and $\mathfrak{A}_2$ can be derived
by \cite[Lemma 3]{Sver93}. Actually, our proof is a simple adaptation.
\end{rem}
\begin{rem}
Since the matrices in $\mathfrak{A}^{co}_2-\dig{(3,1)}$ are not positive definite, by a straightforward
modification of Lemma \ref{sverak split} we infer that the homogenized set related to 
$\lbrace(\mathfrak{A}_1,P_1),(\mathfrak{A}^{co}_2,P_2)\rbrace$ is not polyconvex.
Consequently, the homogenized integrand related to 
$\chi_{P_1}(y)\,\W_1(\Lambda)+\chi_{P_2}(y)\,\mathrm{dist}^p(\Lambda,\mathfrak{A}^{co}_2)$
cannot be polyconvex. This proves that also by mixing a polyconvex function and a convex function both
$p$-coercive and with $p$-growth, loss of polyconvexity can occur in the  homogenization process.
\end{rem}

%%%%%%%%%%%%%%%%%%%%%%%%%%%%%%%%%%%%%%%%%%%%%%%%%%%%%%%%%%%%%%%%%%%%%%%%%%%%%%%%%%%%%%%%%%%%%%%%%%%
%%%%%%%%%%%%%%%%%%%%%%%%%%%%%%%%%%%%%%%%%%%%%%%%%%%%%%%%%%%%%%%%%%%%%%%%%%%%%%%%%%%%%%%%%%%%%%%%%%%

%%%RINGRAZIAMENTI-----------------------------------------------------------------------------------------------------
\bigskip
\centerline{\textsc{\large{Acknowledgments}}}
\bigskip
\noindent
I wish to thank Gianni Dal Maso, Antonio DeSimone and Enzo Nesi for several stimulating discussions on the topic.

%%%BIBLIOGRAFIA--------------------------------------------------------------------------------------------------------

\end{document}